\documentclass{amsart}
\usepackage{color}

\usepackage{amsmath} 
\usepackage{amssymb}
\usepackage{mathrsfs}

\newtheorem{theorem}{Theorem}[section] 
\newtheorem{claim}[theorem]{Claim}
\newtheorem{lemma}[theorem]{Lemma} 
\newtheorem{proposition}[theorem]{Proposition} 
\newtheorem{conclusion}[theorem]{Conclusion}
\newtheorem{obs}[theorem]{Observation}

\theoremstyle{definition}
\newtheorem{definition}[theorem]{Definition}
\newtheorem{observation}[theorem]{Observation}

\newtheorem{fact}[theorem]{Fact}
\newtheorem{discussion}[theorem]{Discussion}

\theoremstyle{remark}
\newtheorem{remark}[theorem]{Remark}
\newtheorem{question}[theorem]{Question}
\newtheorem{notation}[theorem]{Notation}

\newcommand{\st}{{\rm such that}}

\newcommand{\Th}{{\rm Th}}

\newcommand{\add}{{\rm add}}

\newcommand{\BPA}{{\rm BPA}}

\newcommand{\Max}{{\rm Max}}

\newcommand{\mlt}{{\rm mlt}}

\newcommand{\cof}{{\rm cof}}

\newcommand{\id}{{\rm id}}

\newcommand{\PA}{{\rm PA}}

\newcommand{\cat}{{\rm cat}}

\newcommand{\mdt}{{\rm mdt}}

\newcommand{\Min}{{\rm Min}}

\newcommand{\rest}{{\restriction}}

\newcommand{\wilog}{{\rm without loss of generality}}

\newcommand{\then}{{\underline{then}}}
\newcommand{\when}{{\underline{when}}}
\newcommand{\Then}{{\underline{Then}}}

\newcommand{\Iff}{{\underline{iff}}}
\newcommand{\mn}{{\medskip\noindent}}
\newcommand{\sn}{{\smallskip\noindent}}

\newcommand{\cC}{{\mathscr C}}

\newcommand{\cO}{{\mathscr O}}

\newcommand{\bbL}{{\mathbb L}}
\newcommand{\bbN}{{\mathbb N}}

\newcommand{\bbQ}{{\mathbb Q}}
\newcommand{\bbZ}{{\mathbb Z}}

\newcommand{\bbR}{{\mathbb R}}

\newcommand{\cf}{{\rm cf}}

\newcommand{\cut}{{\rm cut}}

\newcommand{\seq}{{sequence}}
\newcommand{\cont}{{continuous}}
\newcommand{\incr}{{increasing}}

\newcommand{\isom}{{isomorphism}}

\newcommand{\dec}{{decreasing}}
\newcommand{\contr}{{contradiction}}

\newcount\skewfactor
\def\mathunderaccent#1#2 {\let\theaccent#1\skewfactor#2
\mathpalette\putaccentunder}
\def\putaccentunder#1#2{\oalign{$#1#2$\crcr\hidewidth
\vbox to.2ex{\hbox{$#1\skew\skewfactor\theaccent{}$}\vss}\hidewidth}}

\newenvironment{PROOF}[2][\proofname.]
   {\begin{proof}[#1]\renewcommand{\qedsymbol}{\textsquare$_{\rm #2}$}}
   {\end{proof}}

\newcount\skewfactor
\def\mathunderaccent#1#2 {\let\theaccent#1\skewfactor#2
\mathpalette\putaccentunder}
\def\putaccentunder#1#2{\oalign{$#1#2$\crcr\hidewidth
\vbox to.2ex{\hbox{$#1\skew\skewfactor\theaccent{}$}\vss}\hidewidth}}

\def\smallbox#1{\leavevmode\thinspace\hbox{\vrule\vtop{\vbox
   {\hrule\kern1pt\hbox{\vphantom{\tt/}\thinspace{\tt#1}\thinspace}}
   \kern1pt\hrule}\vrule}\thinspace}

\usepackage[hidelinks]{hyperref}

\begin{document}
\makeatletter\def\shfiuwefootnote{\gdef\@thefnmark{}\@footnotetext}\makeatother\shfiuwefootnote{Version 2016-04-13\_11. See \url{https://shelah.logic.at/papers/757/} for possible updates.}

\title{Quite Complete Real Closed fields \\
Sh757}
\author{Saharon Shelah}
\address{Institute of Mathematics
 The Hebrew University of Jerusalem
 Jerusalem 91904, Israel
 and  Department of Mathematics
 Rutgers University
 New Brunswick, NJ 08854, USA}
\email{shelah@math.huji.ac.il}
% \urladdr{http://www.math.rutgers.edu/\char`\~shelah}
%%	Notimplemented locally
\thanks{First typed: March 2000. \\
Research supported by the United States-Israel Binational
Science Foundation. Publication 757}

\subjclass[2010]{Primary: 03C64, 03C60; Secondary: 03C55, 03C98, 13L05}

% Previous version: April 11, 2016

\date{April 13, 2016}

\setcounter{section}{0}

\keywords{Real closed fields, ordered fields, cuts,
completeness, cofinality, symmetric cuts, symmetric closure}

\let\labeloriginal\label
\let\reforiginal\ref

\begin{abstract}
We prove that any ordered field can be
extended to one for which every decreasing sequence
of bounded closed intervals, of any length,
has a nonempty intersection; equivalently, there are no Dedekind cuts
with equal cofinality from both sides.
Here we strengthen the results from the published version.
\end{abstract}

\maketitle
\numberwithin{equation}{section}

\newpage

\section{Introduction}

Laszlo Csirmaz raised the question of the existence of non-archimedean
ordered fields with the following completeness property: any
decreasing sequence of closed bounded intervals, of any ordinal length,
has nonempty intersection. We will refer to such fields as {\it
symmetrically complete} for reasons indicated below.

\begin{theorem}
\label{y4}
1) Let $K$ be an arbitrary ordered field.
\Then \, there is a symmetrically complete real closed field $K^+$
containing $K$ such that any asymmetric cut of $K$ is not filled so if
$K$ is not embeddable into $\bbR$ then $K^+,K$ necessarily have 
an asymmetric cut.

\noindent
2) Moreover it is embeddable over $K$ into $K'$ for any symmetrically
   closed $K' \supseteq K$ and is unique up to isomorphisms over $K$.
\end{theorem}

The construction shows that there is even a ``symmetric-closure'' in a
natural sense, and that the cardinality may be taken to be at most
$2^{|K|^+ +\aleph_1}$.

I thank the referee for rewriting the paper as appeared in the Israel Journal.

In September 2005 (after the paper appeared),
lecturing on it in the Rutgers logic
seminar without details, Cherlin asked where the bound
$\kappa\geq d (K)$ for the number of steps needed in the construction
was the true one.
Checking the proof, it appears that this was used
(in the published version) and eventually we show that this is the
right bound.

Note that by \cite{Sh:405}, consistently with ZFC 
(i.e., after forcing extension, in fact just adding enough Cohen
reals) for some non-principal ultrafilter $D$ on
$\bbN,\bbR^{\bbN}/D$, which is (an $\aleph_1$-saturated) ultrapower 
of the field of the reals (hence a real closed field), is Scott complete.
Also compared to the published version we expand \S5 dealing with
other related closures; note also that if $K$ is an
order field which is symmetrically closed or just
has no cut of cofinality $(\cf(K),\cf(K))$ \then\ $K$ is real closed.

Note also that being symmetrically complete is dual to being quite
saturated because if $K$ (a real closed field or just linear order)
which is $\kappa$-saturated and have a $(\theta_1,\theta_2)$-cut then
$\theta_1 < \kappa \Rightarrow \theta_2 \ge \kappa$.

Our problem of constructing such fields translate to considering cuts of $K$
and their pair of cofinalities. Our strategy is:
\mn
\begin{enumerate}
\item[$(a)$]   we consider some properties of cuts (of a real 
closed field), namely being Dedekind, being Scott, being positive, 
being additive, being multiplicative,
\sn
\item[$(b)$] we define dependence relation on the set of cuts
of $K$, which satisfies the Steinitz assumptions,
\sn
\item[$(c)$] realizing a maximal independent family of cuts
with the right pairs of cofinality, we get a ``one step symmetric closure''.
\end{enumerate}
\mn
It is fine: we can show existence and uniqueness.
But will iterating this ``atomic closure'' eventually terminate?
\mn
\begin{enumerate}
\item[$(d)$] For a field $K$ we define a similar chain
inside $K$, its minimal length being called $h(K)$
\sn
\item[$(e)$]  we define the depth $d(K)$ of $K$
\sn
\item[$(f)$] we show that $h(K)$, and $d(K)$ are quite
closed and show that our iterated closure from (c) does not increase 
the cardinality too much
\sn
\item[$(g)$] we finally show that after $d(K)+1$
steps the iterated closure from (c) terminate.
\end{enumerate}
\mn
Two later works (around 2012) seem relevant. One is by
Kuhlmann-Kuhlman-Shelah \cite{Sh:1024} deal with symmetrically
complete ordered sets for generalizations of Banach fix point theorem.

The other, Malliaris-Shelah \cite{Sh:998} prove that for $T = \Th(\bbN)$
this is false.  Malliaris ask (2013) whether we can generalize Theorem
\ref{y4} to any o-minimal theory $T$.  This is very reasonable but by
the proof of \cite[\S2]{Sh:405} for, e.g. the theory of 
$(\bbR,e^x)$ this fails, but it
is o-minimal by the celebrated theorem of Wilkie \cite{Wil96}.  
\bigskip

\noindent
More fully
\begin{theorem}
\label{y6}
1) Any model of $T$ has a symmetric cut \when \,:
\mn
\begin{enumerate}
\item[$(a)$]  let $T$ be a first order theory extending the theory of
ordered semi-rings (so have 0,1 order, addition and multipliation
with the usual rules but $x-y$ does not necessarily exist), which may
have additional symbols (e.g. $e^x$)
\sn
\item[$(b)$]  $T$ implies some first order formula $\varphi(x,y)$
define a function $f_\varphi$ such that:
\sn
\item[${{}}$]  $(\alpha) \quad$ for $x > 0,f_\varphi(x)$ increases and
is $>x$
\sn
\item[${{}}$]  $(\beta) \quad$ for $x > 1$ we have
$xf_\varphi(x-1) < f_\varphi(x)$.
\end{enumerate}
\mn
2) We can weaken (b) to:
\mn
\begin{enumerate}
\item[$(b)'$]  for some formula $\varphi(x,y) \in \bbL(\tau_T)$ the
theory $T$ implies:
\sn
\item[${{}}$]  $(\alpha) \quad 0 < x \rightarrow (\exists
y)\varphi(x,y)$
\sn
\item[${{}}$]  $(\beta) \quad \varphi(x,y) \rightarrow x < y$
\sn
\item[${{}}$]  $(\gamma) \quad 1 < x \wedge \varphi(x,y_1) \wedge
\varphi(x_1,y_1) \rightarrow xy_1 < y_2$.
\end{enumerate}
\end{theorem}

\begin{PROOF}{\ref{y6}}
By the proof of \cite[Th.2.2,pg.377-8]{Sh:405} noticing:
\mn
\begin{enumerate}
\item[$(*)$]  the proof was written for $T=\PA$, Peano Arithmetic, but
we use only:
\sn
\item[${{}}$]  $(a) \quad$ any $M \models T$ is an ordered semi-ring
\sn
\item[${{}}$]  $(b) \quad$ a function called exp which in
\cite{Sh:405} is $x^x$, but only the properties

\hskip25pt  listed in the theorem are used.
\end{enumerate}
\mn
Alternatively it follows from part (2).

\noindent
2) Similarly, but see details in \ref{f2}.
\end{PROOF}

\begin{conclusion}
\label{y8}
The theory $T = \Th(\bbR,e^x)$ is o-minimal (by \cite{Wil96} and) any
model of $T$ has a symmetric cut.
\end{conclusion}

\begin{PROOF}{\ref{y8}}
The function $x \mapsto e^{e^x}$ satisfies the requirement from \ref{y6}.
\end{PROOF}

\noindent
Still
\begin{claim}
\label{y11}
Let $T$ be an o-minimal (complete first order) theory, for notational
transparency with elimination of quantifiers, and let $K,L$ denote
models of it.

\noindent
1) All the Definitions and Claims not mentioning
   Scott/additive/multiplicative cuts hold, that is, we have
   \ref{monotone}, \ref{2.2}(1)-(4), \ref{2.2.1}, \ref{2.3.13} and
   \ref{b24} - \ref{b36} and \ref{c2} - \ref{3.6}.

\noindent
2) So in particular the symmetric hull of a model $K$ of $T$
 \ref{hull} and $\alpha(K)$ (see \ref{b36}) are well defined.

\noindent
3) But \ref{y8} says that $\alpha(T) = \infty$ for $T =
   \Th(\bbR,e^x)$.

\noindent
4) If $\lambda > |K|$ is a Ramsey cardinal or just $\lambda
   \rightarrow (\omega_1)^{< \omega}_{|K|}$ and $\alpha(K) \ge
   \lambda$ \then \, for some $K$ of cardinality $|T|,\alpha(K) = \infty$.
\end{claim}

\begin{discussion}
%\label{}
Note we conjecture that $\lambda$ can be chosen as $\beth_{|K|^+}$ or
so in \ref{y11}(4).
\end{discussion}
\newpage

\section{Real Closed Fields}

Any ordered field embeds into a real closed field, and in fact
has a unique real closure. We will find it convenient to work
mainly with real closed fields throughout.
Accordingly, we will need various properties of real closed fields.
We assume some familiarity with quantifier elimination, real closure,
and the like, and we use the following consequence of $o$-minimality.
(Readers unfamiliar with $o$-minimality in general may simply remain
in the context of real closed fields, or in geometrical language,
{\it semialgebraic geometry}.)

\begin{fact}
\label{monotone}
Let $K$ be a real closed field, and
let $f$ be a parametrically definable function of one variable
defined over $K$. Then $f$ is piecewise monotonic, with each
piece either constant or strictly monotonic;
that is we can find finitely many intervals
(allowing $\infty$ and $-\infty$ as end points),
$f$ is constant or strictly increasing or strictly \dec\ on each interval.
Moreover this holds uniformly and
definably in definable families, with
a bound on the number of pieces required, and with each piece an interval whose
endpoints are definable from the defining parameters for the function.
\end{fact}

\begin{notation}
\label{2.1.1}
1) $K,L$ are ordered fields, usually real closed.

\noindent
2) If $K\subseteq L, X\subseteq L$, then $K(X)$ is the (ordered) subfield
of $L$ generated by $K\cup X$.  

\noindent
3) Let $K_+ = \{c\in K:c>0\}$.

\noindent
4) For $K$ an ordered field and $A,B \subseteq K$ let $A \subseteq B$
   mean $a \in A \wedge b \in B \rightarrow a <_K b$.
\end{notation}

\subsection {Cuts}\
\bigskip

\begin{definition}
\label{2.2}
1) A cut in a real closed field $K$ is a pair $C=(C^-,C^+)$
with $K$ the disjoint union of $C^-$ and $C^+$, and $C^-<C^+$.
The cut is a {\em Dedekind cut} if both sides are nonempty,
and $C^-$ has no maximum, while $C^+$ has no minimum.
For $L \subseteq K$ let $C\rest L= (C^-\cap L, C^+\cap L)$,
it is a cut of $L$.

\noindent
2)  The {\it cofinality} of a cut $C$ is the
pair $(\kappa,\lambda)$ with $\kappa$ the cofinality of $C^-$ and
$\lambda$ the coinitiality of $C^+$ (i.e., the ``cofinality to the left'').
If the cut is not necessarily a Dedekind cut, then one includes $0$ and $1$ as
possible values for these invariants.

\noindent
3)  A cut of cofinality $(\kappa,\lambda)$ is {\em symmetric}
if $\kappa=\lambda$.

\noindent
4) A real closed field is {\em symmetrically complete}
if it has no symmetric cuts.

\noindent
5)  A cut is {\em positive} if $C^- \cap K_+$ is nonempty.
\end{definition}

\begin{conclusion}
\label{2.2.1}
1)  If $L$ is a real closed field extending $K$ and $a,b \in L
\setminus K$ realizes the same cut of $K$ \then\ $a,b$ realizes 
the same types over $K$ in $L$.

\noindent
2)  If $C$ is a non-Dedekind cut of $K$ and $\theta = \cf(K)$
 \then \, $\cf(C) \in
\{(\theta,0),(0,\theta),(\theta,1),(1,\theta)\}$ and each of those
 cases occurs.
\end{conclusion}

\begin{proof}
%{\ref{2.2.1}}
By Fact \ref{monotone}.
\end{proof}

\noindent
We will need to consider some more specialized properties of cuts.

\begin{definition}
\label{2.3A}
Let $K$ be a real closed field, $C$ a cut in $K$.

\noindent
1) The cut $C$ is a {\em Scott cut} if it is a Dedekind cut,
and for all $r>0$ in $K$, there
are elements $a\in C^-$, $b\in C^+$ with $b-a<r$.

\noindent
2)  The cut $C$ is {\em additive} if $C^-$ is closed under addition and
contains some positive element.

\noindent
3)  The cut $C$ is {\em multiplicative} if $C^-\cap K_+$ is closed under
multiplication and contains $2$.

\noindent
4)  $C_\add$ is the cut with left side $\{r\in K:r+C^- \subseteq C^-\}$.

\noindent
5) For $C$ a positive cut, $C_\mlt$ is the cut with left side
$\{r\in K:r\cdot (C^- \cap K_+) \subseteq C^-\}$.
\end{definition}

\noindent
Observe that
\begin{obs}
\label{2.3.1}
1)  Scott cuts are symmetric, in fact both cofinalities are
equal to $\cf(K)$.

\noindent
2)  If $C$ is a Dedekind cut which is not a Scott cut, 
\then \, $C_\add$ is a positive additive cut and 
note: $C^-_{\rm add}=\{c:c\leq 0\}$
is impossible as ``not scott''.

\noindent
3) If $C$ is an additive cut which is not a multiplicative cut, 
\then \, $C_\mlt$ is a multiplicative cut.

\noindent
4) If $C$ is a Dedekind cut of cofinality $(\kappa,\lambda)$ then 
$\kappa,\lambda \ge \aleph_0$.
\end{obs}

\begin{definition}
\label{2.3.13}
1) If $K \subseteq L$ are ordered fields, then a cut $C$ in $K$ is said to
be {\it realized}, or {\it filled}, by an element $a$ of $L$
\Iff \, the cut induced by $a$ on $K$ is the cut $C$.

\noindent
2)  If $C_1,C_2 \subseteq K$ and $C_1<C_2$ but no
$a\in K$ satisfies $C_1<a<C_2$ \then \, the cut of $K$ defined 
(or induce or canonically extends) by $(C_1,C_2)$ is $(\{a\in K:a \le c$ 
for some $c\in C_1\}, \{b\in K:c\leq b$ for some $c\in C_2\})$,
e.g., $(C_1,C_2)$ may be a cut of a subfield of $K$.

\noindent
3) If $K \subseteq L$ and $C$ is a cut of $L$ \then \, $C \rest K =
   (C^- \cap K,C^+ \rest K)$, a cut of $K$, is called the cut of $K$
   induced by $C$.
\end{definition}

\noindent
By Scott \cite{Sco69} we know that
\begin{lemma}
\label{scott}
Let $K$ be a real closed field.  \Then \, there is a real closed field $L$
extending $K$ in which every Scott cut has a unique realization,
and no other Dedekind cuts are filled.
\end{lemma}

This is called the {\it Scott completion} of $K$, and is strictly
analogous to the classical Dedekind completion. The statement
found in \cite{Sco69} is worded differently, without
referring directly to cuts, though
the relevant cuts are introduced in the course of the proof.
The result is also given in greater generality there.

\begin{lemma}
\label{b18}
Let $K$ be a real closed field, $C$ a multiplicative cut in $K$,
and $L$ the real closure of $K(x)$, where $x$ realizes the cut $C$.
\Then\ for any $y\in L$ realizing the same cut, we have
$x^{1/n}<y<x^n$ for some $n$.
\end{lemma}

\begin{PROOF}{\ref{b18}}
%{\ref{multiplicative}}
Let $\cO_K$ be $\{a\in K:|a|\in C^-\}$, and let $\cO_L$ be the convex
closure in $L$ of $\cO_K$. Then these are valuation rings,
corresponding to valuations on $K$ and $L$ which will be called
$v_K$ and $v_L$ respectively.

The value group $\Gamma_K$ of $v_K$ is a divisible ordered abelian
group, and the value group of the restriction of $v_L$ to $K(x)$
is $\Gamma_K \oplus \bbZ \gamma$ where $\gamma := v_L(x)$ is negative, and
infinitesimal relative to $\Gamma_K$. The value group of $v_L$ is the
divisible hull of $\Gamma_K\oplus \bbZ \gamma$.

Now if $y\in L$ induces the same cut $C$ on $K$, then $v_L(y)=qv_L(x)$
for some positive rational $q$. Hence $u=y/x^q$ is a unit of $\cO_L$,
and thus $u,u^{-1}<x^\epsilon$ for all positive rational $\epsilon$.
So $x^{q-\epsilon}<y<x^{q+\epsilon}$ and the claim follows.
\end{PROOF}

\begin{lemma}
\label{additive}
Let $K \subseteq L$ be real closed fields, and $C$ an additive
cut in $L$.
Let $C'$ and $C'_\mlt$ be the cuts induced on $K$ by $C$ and $C_\mlt$
respectively. Suppose that $C'_\mlt=(C')_\mlt$, and
that $x,y\in L$ are two realizations of the cut $C'$, with
$x\in C^-$ and $y\in C^+$. Then $y/x$ induces the cut $C'_\mlt$ on $K$.
\end{lemma}

\begin{PROOF}{\ref{additive}}
%{\ref{additive}}
If $a \in K_+$ and $ax\ge y$, then $a\in (C_\mlt)^+$, by definition, working
in $L$.

On the other hand if $a \in K_+$ and $ax< y$, then
$a \in [(C')_\mlt]^- \cap \kappa$, which by hypothesis is $(C'_\mlt)^-$.
\end{PROOF}

\begin{lemma}
\label{dedekind}
Let $K \subseteq L$ be real closed fields, and $C$ a positive
Dedekind cut in $L$ which is not additive.
Let $C'$ and $C'_\add$ be the cuts induced on $K$ by $C$ and $C_\add$
respectively. Suppose that $C'_\add=(C')_\add$.
Suppose that $x,y\in L$ are two realizations of the cut $C'$, with
$x\in C^-$ and $y\in C^+$. Then $y-x$ induces the cut $C'_\add$ on $K$.
\end{lemma}

\begin{PROOF}{\ref{dedekind}}
%{\ref{dedekind}}
If $a\in K$ and $a+x\ge y$, then $a\in (C_\add)^+$, by definition, working
in $L$.

On the other hand if $a\in K$ and $a+x< y$, then
$a\in [(C')_\add]^- \cap K$, which by hypothesis is $(C'_\add)^-$.
\end{PROOF}
\bigskip

\subsection {Independent cuts}\
\bigskip

We will rely heavily on the following notion of independence.

\begin{definition}
\label{b24}
Let $K$ be a real closed field, and $\cC$ a set of cuts in $K$.
We say that the cuts in $\cC$ are {\em dependent} if for every
real closed field $L$ containing realizations $a_C$ ($C \in \cC$)
of the cuts over $K$, the set $\{a_C:C\in \cC\}$
is algebraically dependent over $K$.
\end{definition}

The following merely rephrases the definition (recalling \ref{2.2.1}(1))
\begin{lemma}
\label{2.9}
Let $K$ be a real closed field and $\cC$ a set of cuts over $K$.

\noindent
1) The following are equivalent:
\mn
\begin{enumerate}
\item[$(A)$]   $\cC$ is independent
\sn
\item[$(B)$]  For each set $\cC_0 \subseteq \cC$, and each ordered field $L$
containing $K$, if $a_C\in L$ is a realization of the cut $C$ for
each $C \in \cC_0$, \then \, the real closure of $K(a_C:C \in \cC_0)$ does not
realize any cuts in $\cC \setminus \cC_0$.
\end{enumerate}
\mn
2)  For ${\cC}_0$ and $L$ as in clause (B) above, every $C \in {\cC}
\setminus {\cC}_0$ define a (unique) cut $C'$ of $L$
(see Definition \ref{2.2}(2)) and $\{C':C\in {\cC} \setminus
 {\cC}_0\}$ is an independent set of cuts of $L$.

\noindent
3)  Assume $\langle K_\alpha:\alpha \leq \delta\rangle$ is an 
increasing \seq\ of real closed fields, ${\cC}$ a set of cuts of $K_0$, and
$C \in {\cC} \wedge \alpha < \delta \Rightarrow {\cC}$ define a cut
$C^{[\alpha]}$ of $K_\alpha$.
\Then\ each $C \in {\cC}$ define a cut of $K_\delta$ which we
call $C^{[\delta]}$ and if $\{C^{[\alpha]}:C \in \cC\}$ is an independent
set of cuts of $K_\alpha$ for each $\alpha \in \delta$ 
then $\{C^{[\delta]}:C \in {\cC}\}$ is an
independent set of cuts of $K_\delta$.

\noindent
4) This dependence relation satisfies the Steinitz axioms for a
dependence relation. 
\end{lemma}

\noindent
 We will make use of it to realize certain sets of
types in a controlled and canonical way.
\begin{lemma}
\label{b28}
Let $K$ be a real closed field, and $\cC$ a set of
cuts over $K$.

\noindent
1) There is a real closed field $L$ generated over
$K$ (as a real closed field) by a set of realizations of some independent
family of cuts included in $\cC$,
in which all of the cuts $\cC$ are realized. 

\noindent
2) Furthermore, such an extension is unique up to isomorphism over
   $K$.

\noindent
3) Assume $\cC = \bigcup\limits_{i < \alpha} \cC_i$ is
   $\subseteq$-increasing continuous.  There is a sequence $\langle
 K_i:i \le \alpha \rangle$ of real closed fields, $K_0 = K,K_{i+1}$
is gotten as in (1) for $(K_i,\cC'_i)$ where $\cC'_i = \{C':C'$ is a
   cut of $K_i$ which is induced by $C' \rest K_0$ and $C' \rest K_0
   \in \cC_i\}$.
Moreover, this sequence is unique up to isomorphism and $K_\alpha$
   is gotten from the pair $(K,\cC)$ as in part (1).
\end{lemma}

\begin{PROOF}{\ref{b28}}
1) Clearly we must take $L$ to be the real closure of $K(a_C:C\in
\cC_0)$, where $\cC_0$ is some maximal independent subset of $\cC$;
and equally clearly, this works.

It remains to check the uniqueness. This comes down to the following: for
any real closed field $L$ extending $K$, and for any choice of
independent cuts $C_1,\dots,C_n$ in $K$ which are realized by elements
$a_1,\dots,a_n$ of $L$, the real closure of the
field $K(a_1,\dots,a_n)$ is uniquely
determined by the cuts. One proceeds by induction on $n$. The real
closure $\hat K$ of $K(a_n)$ is determined by the cut $C_n$
by \ref{2.2.1}; and as none of the
other cuts are realized in it, they extend canonically to cuts
$C_1',\dots,C_{n-1}'$ over $\hat K$, which are independent over $\hat
K$. At this point induction applies.

\noindent
2),3) Easy, too.
\end{PROOF}

\begin{lemma}
\label{b30}
Let $K$ be a real closed field, and
$\cC$ a set of Dedekind cuts in $K$. Suppose that $C$ is a
Dedekind cut of $K$ of cofinality $(\kappa,\lambda)$ which is dependent on
$\cC$, and let $\cC_0$ be the set $\{C' \in
\cC:\cof(C')=\hbox{$(\kappa,\lambda)$ or $(\lambda,\kappa)$}\}$.
\Then\footnote{If we allow $\{\kappa,\lambda\} \cap \{0,1\} \ne
\emptyset$ then we should equate 0 and 1.} \, $C$ is dependent 
on $\cC_0$, and in particular $\cC_0$ is non-empty.
\end{lemma}

\begin{PROOF}{\ref{b30}}
It is enough to prove this for the case that $\cC$ is
independent. If this fails, we may replace the base 
field $K$ by the real closure
$\hat K$ over $K$ of a set of realizations of $\cC_0$. Then since
none of the cuts in $\cC \setminus \cC_0$ are realized, and $C$ is not
realized, these cuts extend canonically to cuts over $\hat K$, and
hence we may suppose $\cC_0=\emptyset$. We may also suppose $\cC$ is
finite, and after a second extension of $K$ we may even assume
that $\cC$ consists of a single cut $C_0$. This is the essential case.

So at this point we have $L \supseteq K$ and in it a 
realization $a$ of $C_0$ over the real
closed field $K$, and a realization $b$ of $C$ over $K$, with $b$
algebraic, and hence definable, over $a$, relative to $K$.
Thus $b$ is the value at $a$ of a $K$-definable function
 (being the $\ell$-th root for some $\ell$ and polynomial over $K$), 
not locally constant near $a$, and by Fact \ref{monotone} it follows 
that there is an interval $I_0$ of $L$ containing 
$b$ with endpoints in $K$ and a $K$-definable function
which is order isomorphism or anti-isomorphism from
the interval $I_0$ to an interval $I_1$ including $a$, with the 
cuts corresponding. So the cuts have same (or inverted) cofinalities.
This contradicts the supposition that $\cC_0$ has become empty, and 
proves the claim.
\end{PROOF}

\noindent
For our purposes, the following case is the main one. We combine
our previous lemma with the uniqueness statement.

\begin{proposition}
\label{symmetry}
Let $K$ be a real closed field, and $\cC$ a maximal independent set of
symmetric cuts in $K$. Let $L$ be an ordered field containing $K$
together with realizations $a_C$ of each $C \in \cC$.
\Then\ the real closure $K'$
of $K(a_C:C \in \cC)$ realizes the symmetric cuts of $K$ and no others.
Furthermore, the result of this construction is unique up to
isomorphism over $K$.  Moreover if $L'$ is a real closed field extending $K$
which realizes every cut in
$\cC$ {\bf then} $K'$ can be embedded into $L'$ over $K$.
\end{proposition}

\begin{PROOF}{\ref{symmetry}}
Clear (the ``no other'' by \ref{b30}).
\end{PROOF}

\noindent
Evidently, this construction deserves a name.
\begin{definition}
\label{hull}
1)  Let $K$ be a real closed field. A {\em symmetric hull} of $K$ is a
real closed field generated over $K$ by a set of realizations of a
maximal independent set of symmetric cuts.

\noindent
2)  We say that $\bar K=\langle K_\alpha:\alpha\leq\alpha^*\rangle$
is an associated symmetric $\alpha_*$-chain over $K$ when:
\mn
\begin{enumerate}
\item[$(a)$]  $K_0=K$,
\sn
\item[$(b)$]  $K_\alpha$ is an order field,
\sn
\item[$(c)$]  $K_\alpha$ is increasing continuous with $\alpha$
\sn
\item[(d)] $K_{\alpha+1}$ is a symmetric hull of $K_\alpha$
for $\alpha<\alpha^*$.
\end{enumerate}
\mn
3)  In (2) we replace ``$\alpha_*$-chain'' by ``chain'' if 
$K_{\alpha_*}$ is symmetrically complete but 
$\alpha < \alpha_* \Rightarrow K_{\alpha+1} \neq K_\alpha$.

\noindent
4) Let $\alpha(K)$ be the minimal $\alpha_*$ such that for some
   $\langle K_\alpha:\alpha \le \alpha_*\rangle$ as in part (2) is a
   chain, and $\infty$ if there is no such $\alpha$.
\end{definition}

\begin{conclusion}
\label{b36}
1) For every $K$, some $L$ is a symmetric hull of $K$, and it is
   unique up to isomorphisms over $K$.

\noindent
2) For every $K$ and $\alpha_*$ there is an associated symmetric
$\alpha_*$-chain $\langle K_\alpha:\alpha \le \alpha_*\rangle$ over
$K$.  It is unique up to isomorphism over $K$, that is, if
   $\langle K'_\alpha:\alpha \le \alpha_*\rangle$ is another such
   $\alpha^*$-chain, then there is an isomorphism $f$ from
   $K_{\alpha_*}$ onto $K'_{\alpha_*}$ such that $f \rest K = \id_K$
and $f$ maps $K_\alpha$ onto $K'_\alpha$ for $\alpha \le \alpha_*$.

\noindent
3) $\alpha_* = \alpha(K) < \infty$ iff for every $\beta > \alpha_*$
   and associated symmetric $\beta$-chain $\bar K$ over $K$ we have
   $(\forall \gamma)(\alpha_* \le \gamma \le \beta \rightarrow
K_\gamma = K_{\gamma_+ +1}) \wedge (\forall \gamma)(\gamma < \alpha_*
   \rightarrow K_\gamma \ne K_{\gamma +1})$.
\end{conclusion}

\noindent
While the ``symmetric hull" (from \ref{hull}(1)) is unique 
up to isomorphism, there is certainly no reason
to expect it to be symmetrically complete, and the construction will
need to be iterated. The considerations of the next section will
help to prove that the construction eventually terminates and to
bound the length of the iteration.

\begin{lemma}
\label{scottsymmetric}
1)  For regular $\kappa<\lambda$ there is a real closed field
with an $(\kappa,\lambda)$-cut.

\noindent
2)  Let $K$ be a real closed field, and $L$ its symmetric hull.
\Then \, every Scott cut in $K$ has a unique realization in $L$.

\noindent
3)  Assume $L$ is the real closure of $K \cup \{a_C:C \in {\cC}\},{\cC}$
an independent set of cuts of $K$ and $a_C$ realizing
the cut $C$ in $L$ for $C\in {\cC}$
\mn
\begin{enumerate}
\item[$(a)$]  If every $C\in {\cC}$ is a Dedekind cut \then\ every cut 
of $K$ realized in $L$ is a Dedekind cut
\sn
\item[$(b)$]  if  every $C \in \cC$ is a non-Dedekind cut of $L$ 
 \then\ every cut of $K$ is realized in $L$ is non-Dedekind
\sn
\item[$(c)$]  if some $C \in \cC$ is a non-Dedekind cut of $K$ \then
\, every non-Dedekind cut of $K$ is realized in $L$.
\end{enumerate}
\end{lemma}

\begin{PROOF}{\ref{scottsymmetric}}
1)  First we choose $K_i$ a real closed field $K_i$ \incr\ \cont\ with 
$i\leq \kappa,K_0=K$ and for $i<\kappa$ the element
$a_i\in K_{i+1} \setminus K_i$ is above all members of $K_i$.
Second we choose a real closed field $K^i$ \incr\ \cont\
with $i\leq\lambda$ \st\ $K^0=K_\kappa$, and
for $i<\lambda,b_i\in K^{i+1}_i\setminus K_i$
is above $a_j$ for $j<\kappa$ and below
any $b\in K^i$ \st\ $(\forall j<\kappa)
(a_j < b)$, in $K^i$. Lasly in $K^\lambda$,
$(\{a_j:j<\kappa\}, \{b_i:i<\lambda\})$ determine a
$(\kappa,\lambda)$-cut, i.e., ($\{a:a<a_j$ for some $j<\kappa\},
\{b:b_i<b$ for $i<\lambda\}$), in $K_\lambda$, is such a cut.

\noindent
2) Recall that every  Scott cut is symmetric.
One can form the symmetric hull of $K$ by first taking its Scott
completion $K_1$, realizing only the Scott cuts (uniquely), and then
taking the symmetric hull of $K_1$; this is equivalent by
\ref{b28}(3).  By part (3) we are done.

\noindent
3)  Easy by now; Clause (a) is really from \cite{Sco69}.
\end{PROOF}

\begin{observation}
\label{b42}
1) For every linear order $I$ there is a real closed field $L$ and
   order preserving function from $I$ into $K$ such that: for every
   Dedekind cut $(C_1,C_2)$ of $I$, the pair $(\{f(s):s \in
   C_1\},\{f(s):s \in C_2\})$ induce a Dedekind cut of $K$; also
   $|L|=|I|$.

\noindent
2) So, e.g. for every $\mu$ for some $K$ of cardinality $\mu,K$ has a
   $(\theta_1,\theta_2)$-cut whenever $\theta_1,\theta_2 \le \mu$ are regular.
\end{observation}

\begin{PROOF}{\ref{b42}}
1) Let $K$ be any field.  We can find $L \supseteq K$ such that $L =
   K(\{a_s:s \in I\})$ such that $s < t \Rightarrow
   \bigwedge\limits_{n} L \models ``(a_s)^n < a_t"$.  Now $L$ is as
   required and $|L|=|K| + |I|$.

\noindent
2) Easy.
\end{PROOF}
\newpage

\section{Height and Depth}

\begin{definition}
\label{c2}
Let $K$ be a real closed field.

\noindent
1) The {\em height} of $K$, $h(K)$,
is the least ordinal $\alpha$ for which we can find a continuous
increasing sequence $K_i$ ($i\le \alpha$) of real closed fields
with $K_0$ countable, $K_\alpha=K$, and $K_{i+1}$ generated over
$K_i$, as a real closed field, by a set of realizations of a family of
 cuts which is independent.

\noindent
2)  Let $h^+(K)$ be $\max(|h(K)|^+,\aleph_1)$
\end{definition}

\begin{remark}
\label{c5}
1)  $\aleph_1$ is the first uncountable cardinal.

\noindent
2)  $h^+ (K)$ is the first uncountable
cardinal strictly greater than $h(K)$), so regular.

\noindent
3) We could have chosen $K_0$ as the algebraic members of $K$, but
   this is not enough to make $\alpha$ unique.  The point is that
   there may be, e.g. $\langle x_q:q \in \bbQ\rangle$ in $K$ such that
   $q_1 < q_2 \Rightarrow \bigwedge\limits_{n} (x_{q_1})^n < x_{q_2}$,
   so for every $\alpha < \aleph_1$ there is an increasing sequence
   $\langle q_\beta:\beta < \alpha\rangle$ of rationals, so may be
   $x_{q_\beta} \in K_\gamma \leftrightarrow \beta < \gamma$.
   Similarly for any $\lambda > \aleph_0$.

\noindent
4) Observe that the height of $K$ is an ordinal of cardinality
at most $|K|$ (or is undefined, you can let it be
$\infty$, a case which by \ref{3.3} does not occur).
We need to understand the relationship of the height of $K$ and of
$\alpha(K)$ with its order-theoretic structure, which for 
our purposes is controlled by the following parameter.
\end{remark}

\begin{definition}
\label{c8}
Let $K$ be a real closed field.
The {\em depth} of $K$, denoted $d(K)$, is the least cardinal
$\kappa$ greater than the length of every strictly increasing sequence
in $K$.
\end{definition}

\begin{obs}
\label{c11}
If $K$ is a real closed field, \then \, $d(K)$ is a regular uncountable
cardinal.
\end{obs}

\begin{PROOF}{\ref{c11}}
Uncountable because there is an infinite increasing sequence: $1,2,\ldots$.
Regular as any interval of $K$ is order isomorphic to $K$.
\end{PROOF}

\noindent
The following estimate is straightforward, and what we really
need is the estimate in the other direction, which will be given
momentarily.
\begin{lemma}
\label{3.3}
Let $K$ be a real closed field.
\Then \, $h(K)\le d(K)$.
\end{lemma}

\begin{PROOF}{\ref{3.3}}
One builds a continuous strictly \incr\ tower $K_\alpha$ of real 
closed subfields of $K$ starting with any countable subfield of $K$.
If $\alpha$ is limit, we define $K_\alpha = 
\bigcup\limits_{\beta<\alpha} K_\beta$.  For successor ordinals,
$K_{\alpha+1}$, is the real closure of $K_{\alpha} \cup
\{a^\alpha_C:C \in \cC_\alpha\}$ inside $K$ where
${\cC}_\alpha$ is a maxmimal independent set of cuts
of $K_\alpha$ realized in $K$, and
$a^\alpha_C \in K$ realized $C$. We stop when $K_\alpha=K$.  Now if
$K_\alpha \ne K$ then every $a \in K \backslash K_\alpha$ realizes
some cut of $K_\alpha$ so by \ref{2.9}(4), there is $\cC_\alpha$ as
required (and it is not empty) hence $K_{\alpha +1}$ as required can
be chosen.  As $K_\alpha$'s definition implies $|\alpha| \le
|K_\alpha| \le |K|$, necessarily for some $\alpha,K_\alpha = K$ and
then we stop.
If this continues past $\kappa=d(K)$, then there is a cut
over $K_\kappa$ filled at stage $\kappa$ by an element $x\in K$.
Then the cut determined by $x$ over each $K_\alpha$ for
$\alpha<\kappa$ is filled at stage $\alpha+1$ by an element
$y_\alpha$.  Those $y_\alpha$'s lying below $x$ form an
increasing sequence, by construction, which is therefore of length
less than $\kappa$; and similarly (using the $\langle y_\alpha:\alpha
< \kappa\rangle$) there are fewer than $\kappa$
elements $y_\alpha>x$, so we arrive at a contradiction.
\end{PROOF}

\begin{proposition}
\label{3.4}
Let $K$ be a real closed field.
Then $d(K)\le h^+(K)$.
\end{proposition}

\begin{PROOF}{\ref{3.4}}
Let $\kappa>h(K)$ be regular and uncountable, so it suffices to prove
$d(K) \le \kappa$.  Let
$K_\alpha$ ($\alpha\le h(K)$) be a continuous increasing chain of real
closed fields, with $K_0$ countable, $K_{h(K)}=K$, and $K_{i+1}$
generated over $K_i$, as a real closed field, by a set of realizations
of an independent family of cuts.

For $\alpha\le h(K)$ and $X \subseteq K$, let $K_{\alpha,X}$ be the
real closure of $K_\alpha(X)$
inside $K$. We recast our claim as follows to allow
an inductive argument
\mn
\begin{enumerate}
\item[$\circledast$]  For $X \subseteq K$ with $|X|<\kappa$, 
and any $\alpha\le h(K)$, we have $d(K_{\alpha,X})\le \kappa$.
\end{enumerate}
\mn
Now this claim gives the promised result for $\alpha = h(K)$, 
is trivial for $\alpha=0$ as $K_0$ is countable so $K_{\alpha,X}$ has
cardinality $< \kappa$ (for $X \subseteq K,|X| < \kappa$), and
the claim passes smoothly through limit ordinals up to $h(K)$ (because
$\kappa = \cf(\kappa) > h(K))$, so
we need only to consider the passage from $\alpha$ to $\beta=\alpha+1$.
So $K_\beta$ is $K_{\alpha,S}$ with $S$ a set of realizations of an
independent family of cuts over $K_\alpha$, (no two realizing the same
cut, of course), and similarly $K_{\beta,X}$ is
$K_{\alpha,X \cup S}$.

Consider the claim in the following form:

\[
\hbox{$d(K_{\alpha,X \cup S_0})\le\kappa$
for $S_0 \subseteq S$}
\]

\mn
In this form for $S_0 = S$ we get the desired inductive step, and it
clearly holds if $|S_0|<\kappa$, as it is included in the
inductive hypothesis for $\alpha$, and the case $|S_0|\ge \kappa$
reduces at once to the case $|S_0|=\kappa$. So we now assume that
$S_0$ is a set of realizations of an independent family of 
cuts of $K_\alpha$ of cardinality $\kappa$ (one element per cut).

By \ref{2.9}(2),(3),(4) we can find a subset $S_1$ of $S_0$ of cardinality
$\aleph_0+|X|$ \st:
\mn
\begin{enumerate}
\item[$(a)$]  if $s\in S_0\setminus S_1$ then the cut $C$ which
$s$ induce on $K_\alpha$ is not realized in the real closure
$K'_\alpha(\subseteq K)$ of $K_\alpha (X \cup S_1)$
\sn
\item[$(b)$]  the cuts which the $s\in S_0\setminus S_1$ induce on
$K'_\alpha$ form an independent family.
\end{enumerate}
\mn
Let $Y = X \cup S_1$, so $K'_\alpha = K_\alpha(Y),|Y| < \kappa$.

Let $\{s_\epsilon:\epsilon<\kappa\}$ list $S_0 \setminus S_1$.
For $\zeta\leq\kappa$, let $L_\zeta=K_{\alpha, Y \cup
\{s_\epsilon:\epsilon<\zeta\}}$ and let $L=L_\kappa$.
By the induction hypothesis (for $\alpha$ and $X\cup
\{s_\epsilon:\epsilon<\zeta\}$)
we have $d(L_\zeta)\le \kappa$ for $\zeta<\kappa$,
and we shall prove $d(L)\le \kappa$ thus finishing the proof.

Let $C_i$ be the cut realized by $s_i$ over $L_0$. Note that $C_i$
extends canonically to a cut $C^j_i$ on $K_j$ for all $j\le i$, and
for fixed $j$, the set
$\{C^j_i:i\in [j,\kappa)\}$ of cuts of $K_j$ is independent.

Now suppose, toward a contradiction, that there is $\langle
a_i:i<\kappa \rangle$ an increasing sequence in $L$.
For each $i<\kappa$ let the ordinal 
$f(i) < \kappa$ be minimal such that $a_i \in L_{f(i)}$, so necessarily
$f:\kappa\rightarrow \kappa$ is well defined
and for each $j<\kappa$ the set $\{i<\kappa:a_i \in L_j$ equivalently
$f(i) \le j\}$ is a 
bounded subset of $\kappa$ (because $d(L_j) \le \kappa$).
Now for $\varepsilon,i < \kappa$ with $f(i)>\epsilon$
let $B_i^{\epsilon}$ denote the cut induced on
$L_\epsilon$ by $a_i$.  For $i_1 < i_2 < n$ such that $\varepsilon <
f(i_1),f(i_2)$ clearly $a_{i_1} <_K a_{i_2}$ hence
$(B^\varepsilon_{i_1})^- \subseteq (B^\varepsilon_{i_2})^-$. 
With $\epsilon$ held fixed,  and with $i$ varying, as
$d(L_\epsilon)\le\kappa$ we find that the cuts $B_i^\epsilon$
stabilize for large enough $i<\kappa$
(and furthermore, $a_i\notin L_\epsilon$).
Accordingly, for each $\epsilon$ we may select
$j_\epsilon<\kappa$ above $\epsilon$ such
that the cuts $B^\epsilon_i$ coincide for all $i\ge j_\epsilon$.

Now fix a limit ordinal $\delta<\kappa$ such that for all $\epsilon<\delta$ we
have $j_\epsilon<\delta$. We may also require that $a_i\in L_\delta$ for
$i<\delta$. 

Let $\zeta < \kappa$ be such that $a_\zeta \notin K_\delta$, it is well
defined as $d(K_\delta) \le \kappa$ and is $\ge \delta$ as $i < \delta
\Rightarrow a_i \in L_\delta$.

Now $a_\zeta$ is algebraic over $L_\delta(s_i:i\in I_0)$
for some finite subset $I_0$ of $[\delta,\kappa)$ and $I_0 \ne
\emptyset$ because $a_\zeta \notin L_\delta$.  Hence
$a_\zeta$ is algebraic also over $L_\epsilon(s_i:i\in I_0)$ for
some $\epsilon<\delta$. Thus the cut $B^\epsilon_\zeta$ depends on the
cuts $C^\epsilon_i$ ($i\in I_0$) over $L_\epsilon$.
As $j_\epsilon<\delta < \zeta$ necessarily
$B^\epsilon_\zeta = B^\epsilon_{j_\epsilon}$ is realized in $L_\delta$ and
it follows that this cut is also dependent on the sets
$\{C^\epsilon_i:i<\delta$ and $i \ge \varepsilon\}$ of cuts over $L_\epsilon$.
But the cuts $C^\epsilon_i$ for $i\ge \epsilon$ are supposed to 
be independent over $L_\epsilon$, a contradiction.
\end{PROOF}

\begin{proposition}
\label{3.5}
Let $K$ be a real closed field. \Then\
$|h(K)|\le |K|\le 2^{|h(K)|}$.
\end{proposition}

\begin{PROOF}{\ref{3.5}}
The first inequality is clear. For the second, let $\alpha=h(K)$,
$\kappa=|\alpha|+\aleph_0$,  and
let $K_i$ $(i\leq\alpha)$ be a chain of the sort afforded by the definition
of the height. Note that $h^+ (K)=\kappa^+$.
We show by induction on $i$ that $|K_i|\le 2^\kappa$.
Only successor ordinals $i=j+1$ require consideration, where we
suppose $|K_j|\le 2^\kappa$.

Each generator $a$ of $K_i$ over $K_j$
corresponds to a cut $C_a$ in
$K_j$, and each such cut is determined by the choice of some
cofinal sequence $S_a$ in $C_a^-$. Such a sequence $S_a$ may be taken to
have order type a regular cardinal, and will have length less
than $d(K)$. Since $d(K)\leq h^+(K)$ by \ref{3.4}, we 
find that the order type of $S_a$
is at most $\kappa$. So the number of such sequences is at most
$\sum_{\lambda\le \kappa}|K_j|^\lambda\le \kappa
\times (2^\kappa)^\kappa = 2^\kappa$.
\end{PROOF}

\begin{claim}
\label{3.6}
1) Assume $L$ is the Scott completion of $K$
(i.e., in our terms as in definition \ref{2.3A}(1)).
Let ${\cC}$ be a maximal set of Scott cuts of $K$ which is
independent, \then \, $K$ is dense in $L$.

\noindent
1A) Moreover for every Dedekind cut $C$ of $L$ the cut
$C\rest K$ induce $C$ hence $\cf(C)=\cf (C\rest K)$.

\noindent
1B) Moreover every non-Scott cut $C_1$ of $K$ induce a non-Scott cut
$C_2$ of $L$ so $C_1 = C_1 \rest K,\cf(C_2) = \cf(C_1)$.

\noindent
2)  If $K$ has no symmetric cut except possibly
Scott cuts, and $L$ is a Scott completion of $K$ \then \, $L$
is symmetrically complete hence $L$ is a symmetric closure of $K$.
\end{claim}

\begin{PROOF}{\ref{3.6}}
Part (1) by \cite{Sco69}; the others clear, too.
\end{PROOF}
\newpage

\section{Proof of the Theorem}

We now consider the following construction. Given a real closed field
$K$, we form a continuous increasing chain $K_\alpha$ by setting
$K_0=K$, taking $K_{\alpha+1}$ to be the symmetric hull of $K_{\alpha}$ in the
sense of Definition \ref{hull}, and taking unions at limit ordinals.

If at some stage $K_\alpha$ is symmetrically complete, that is
$K_\alpha=K_{\alpha+1}$, then we have the desired symmetrically
complete extension of $K$, and furthermore our extension is prime in a
natural sense. We claim in fact:
\begin{proposition}
\label{termination}
1)  Assume $K$ is a real closed field, $\kappa = h^+(K) + \aleph_2$
and $K_\alpha$ ($\alpha\le\kappa+1$) is an associated continuous 
symmetric $\kappa$-chain, \then\:
\mn
\begin{enumerate}
\item[$(a)$]  if $\cf(K)\neq \kappa$ then $K_\kappa$ is symmetrically complete
so $K_{\kappa+1}=K_\kappa$
\sn
\item[$(b)$] $K_{\kappa+1}$ is symmetrically closed
\sn
\item[$(c)$]  if $\cf(K) = \kappa$ then every Dedekind cut of $K_\kappa$ is a
Scott cut.
\end{enumerate}
\mn
2)  Also
\mn
\begin{enumerate}
\item[$(i)$]  $|K_\kappa| \le 2^{h^+(K)+\aleph_1}$, and
$|K_{\kappa+1}| \le 2^{h^+(K)+\aleph_1}$
\sn
\item[$(ii)$]  if $K'$ is a symmetrically complete extension
of $K$ \then \, $K_{\kappa+1}$ can be embedded into $K'$ over $K$
\sn
\item[$(iii)$]  $K$ is unbounded in $K_{\kappa+1}$ and no non-symmetric
Dedekind cut of $K$ is realized in $K_\kappa$
\sn
\item[$(iv)$]  any two real closed fields extending $K$ which
are symmetrically complete and embeddable into $K_{\kappa +1}$,
are \isom\ over $K$ (so we can say $K_{\kappa +1}$ is the closure)
\sn
\item[$(v)$]  for some unique $\alpha^* \le \kappa +1$
there is an associated continuous chain over $K$ of length $\alpha^*$.
\end{enumerate}
\end{proposition}

The proof of Proposition \ref{termination}
occupies the remainder of this section.

\begin{lemma}
\label{induced}
Suppose that $K$ is a real closed field, and that
$\langle K_\alpha:\alpha\leq \alpha(*)\rangle$
is a continuous chain of iterated
symmetric hulls of any length.
Let $x\in K_{\alpha} \setminus K_\gamma$ with $\alpha > \gamma \ge 0$ 
arbitrary.  \Then \, the cut induced on $K_\gamma$ by $x$ is symmetric.
\end{lemma}

\begin{remark}
\label{c7}
If we use $\kappa = \max\{h{}^+(K),\aleph_2\}$ then
 $\kappa\ge \aleph_2$ is  regular and greater than $h(K)$; in
particular $\kappa\ge d(K)$ by \ref{3.4}.
Furthermore, as $\kappa>h(K)$, we can view the chain $K_\alpha$ as
a continuation of a chain $\hat K_i$ ($i\le h(K)$) of the sort occurring
in the definition of $h(K)$, with $\hat K_{h(K)}=K_0$;
then the concatenated chain gives a
construction of $K_\alpha$ of length at most $h(K)+\alpha<\kappa$, and
hence $h(K_\alpha)<\kappa$ for all $\alpha<\kappa$, and in particular
$d(K_\alpha)\le \kappa$ for all $\alpha<\kappa$ by \ref{3.4}.
\end{remark}

\begin{PROOF}{\ref{induced}}
Let $\beta < \alpha$ be minimal such that the cut in question is filled in
$K_{\beta+1}$. Then the cut induced on $K_{\beta}$ by $x$ is the
canonical extension of the cut induced on $K_\gamma$ by $x$,
and is symmetric by Proposition \ref{symmetry}.
\end{PROOF}

\noindent
We now begin the proof by contradiction of Proposition \ref{termination}(1).
First, we assume (this does not contradict \ref{termination}(1)) that:
\mn
\begin{enumerate}
\item[$\boxplus_1$]  $(a) \quad$ the chain $\bar K = \langle K_\alpha:
\alpha \le \kappa +1\rangle$ is strictly increasing at every step

\hskip25pt  up to $K_\kappa$, and
\sn
\item[${{}}$]  $(b) \quad C$ is a Dedekind cut of $K_\kappa$.
\end{enumerate}
\mn
Now let
\mn
\begin{enumerate}
\item[$\boxplus_2$]  for $\alpha < \kappa$, let 
$C_\alpha$ denote the cut induced on $K_\alpha$ by $C$.
\end{enumerate}
\bigskip

\begin{lemma}
\label{cutsymmetry}
Assume (in addition to $\boxplus_1$) that $C_\alpha$ does not define
$C$ for $\alpha < \kappa$.

\noindent
1)  For any $\alpha<\kappa$, the cut $C_\alpha$ is symmetric, in
particular, a Dedekind cut.

\noindent
2) For every $\alpha<\kappa, C^-_\alpha$ is bounded in $C^-$ 
and $C^+_\alpha$ is bounded in $C^+$ from below.
\end{lemma}

\begin{PROOF}{\ref{cutsymmetry}}
1) Suppose $C_\alpha$ is not symmetric.  Then the cut $C_\alpha$ 
is not realized in $K_\kappa$, by Lemma \ref{induced}.
Hence the cut $C$ is the canonical extension of $C_\alpha$ to
$K_\kappa$, contradicting the Lemma's assumption.

\noindent
2) Toward contradiction assume $C^-_\alpha$ is unbounded 
from below in $C^-$; so necessarily $\beta\in [\alpha,\kappa) \Rightarrow
\cf(C^-_\beta)=\cf (C^-_\alpha)$.

For every $\beta\in [\alpha,\kappa)$, as $K_{\beta +1}$ is a symmetric 
hull of $K_\beta$, by part (1) it follows that some $a_\beta\in
K_{\beta +1}$ realizes the cut $C_\beta$ hence by the present
assumption toward contradiction, $a_\beta \in C^+$ and $a_\beta \in
C^+_\gamma$ for $\gamma \in (\beta,\kappa)$.  
So for every limit $\delta < \kappa$ which is $>
\alpha,\{a_\beta:\beta\in [\alpha,\delta)\}$ is a subset of
$C^+_\delta$ unbounded from below, hence 
$\cf(C_\delta)=(\cf(C^-_\delta),\cf(C^+_\delta)) = 
(\cf(C^-_\alpha),\cf(\delta))$.
As $\aleph_0 \ne \aleph_1$ are regulars $<\kappa$ for some 
$\delta \in (\alpha,\kappa)$ we have $\cf(C^-_\alpha)\neq \cf(\delta)$ hence
$C_\delta$ is not symmetric, contradicting part (1).  So indeed for
$\alpha <\kappa,C^-_\alpha$ is bounded in $C^-$; similarly
$C^+_\alpha$ is unbounded from below in $C^+$.
\end{PROOF}

\noindent
After these preliminaries, (continue the proof of \ref{termination},
for this) we prove:
\mn
\begin{enumerate}
\item[$\boxplus_2$]  if $C$ is symmetric then $C$ is a Scott cut and
$\cf(K) = \kappa$.
\end{enumerate}
\mn
We divide the analysis of the supposed cut $C$ into a number of cases,
each of which leads to a contradiction or to the desired conclusion.

So assume $C$ is symmetric.
\medskip

\noindent
\underline{Case I}:  $C$ is a Scott cut

If $C$ is (a Scott cut and) $\cf(K) = \kappa$, there is nothing 
to be proved, so assume $\cf(K) \neq \kappa$.

By \ref{cutsymmetry}(2), we can find
$\langle a^-_\alpha,a^+_\alpha:\alpha<\kappa\rangle$,
\st\ $a^-_\alpha\in C^-, a^+_\alpha\in C^+$ both realizing
the cut $C_\alpha$. For some club $E$ of $\kappa$
consisting of limit ordinals we have 
$\alpha<\delta\in E \Rightarrow a^-_\alpha,a^+_\alpha\in K_\delta$.
As $C$ is a Scott cut of $K_\kappa$ by the case assumption necessarily
$\langle a^+_\alpha-a^-_\alpha:\alpha<\kappa\rangle$
is a decreasing \seq\ of positive members of $K_\kappa$ with no 
positive lower bound, so
$\langle 1/(a^+_\alpha-a^-_\alpha):\alpha<\kappa\rangle$
is \incr\ cofinal in $K_{\kappa}$, so $\cf(K_\kappa)=\kappa$.
But $K$ is cofinal in $K_\kappa$ hence $\cf(K)=\kappa$, contradicting 
what we have assumed in the beginning of the case.
\medskip

\noindent
\underline{Case II: $C$ is a multiplicative cut}

Let $\alpha<\kappa$ have uncountable cofinality (recall $\kappa \ge \aleph_2$).

The cut $C_\alpha$ is realized in $K_{\alpha+1}$ by some element $a$.
As $C$ is multiplicative, either all positive
rational powers of $a$ lie in $C^-$, or all positive rational powers
of $a$ lie in $C^+$.

On the other hand, $K_{\alpha+1}$ may be constructed in two stages
as follows. First, let $K_{\alpha +1}$ be the real closure of
$K_\alpha(a_{C'}:C' \in \cC \rangle$ where $a_{C'} \in K_{\alpha +1}$ realizes
$C'$ for $C' \in \cC$ and $\cC$ is an independent set
of symmetric cuts in $K_\alpha$ such that $C_\alpha \in \cC$ and
$a_{C_\alpha} = a$.  Let $\cC' = \cC \backslash \{C\}$ and let
$K'_\alpha$ be the real closure of $K_\alpha(a_{C'}:C' \in \cC')$; 
then take the real
closure of $K_\alpha'(a)$, noting that $a$ fills the canonical extension
of the cut $C_\alpha$ to $K_\alpha'$.  By the choice of $\langle
a_{C'}:C' \in \cC\rangle$ clearly the real closure of $K'_\alpha(a)$
is $K_{\alpha +1}$.  
As seen in Lemma \ref{b18}, there are only two cuts
which may possibly be induced by $C$ on $K_{\alpha+1}$ (which has to
be multiplicative), one has lower part $C^- \cap K'_\alpha$ and the other has upper $C^+ \cap
K'_\alpha$.  Now each of those cuts has countable cofinality from 
one side, and uncountable cofinality from the other.

So $C_{\alpha+1}$ is not symmetric, and this is a contradiction to
\ref{cutsymmetry}. 
\medskip

\noindent
\underline{Case III: $C$ is an additive cut}

By \ref{2.3.1}(3) we know that $C_{\mlt}$ is a multiplicative cut (of
$C_\kappa$). 

Let $\langle b^-_\alpha:\alpha<\kappa\rangle$ be \incr\
cofinal in $C^-\cap K_+$ and let
$\langle b^+_\alpha:\alpha<\kappa\rangle$ be decreasing
unbounded from below in $C^+$.
Clearly $\alpha<\kappa\Rightarrow b^+_\alpha / b^-_\alpha
\in (C_{\rm mlt})^+$ and easily $\langle b^+_\alpha/
b^-_\alpha:\alpha<\kappa\rangle$ is a decreasing \seq\
of members of $(C_{\rm mlt})^+$ unbounded from below in it

According to the property of
$\langle b^+_\alpha / b^-_\alpha:\alpha<\kappa\rangle$
the cofinality of
$(C_\mlt)$ from the right is $\kappa$.

Now if the cofinality of $C_\mlt$ from the left is also $\kappa$,
then by \ref{2.3.1}(3) we
contradict Case II. On the other hand if the cofinality of $C_\mlt$ from
the left is $\theta$ which is less than $\kappa$, then from some
point downward this cofinality stabilizes.  Hence for some closed
unbounded set $E \subseteq \kappa$ we have $\delta \in E \Rightarrow
\cf(C_{\mdt} \rest K_\beta) = (\theta,\cf(\delta))$; but then we can
choose $\delta$ large and of some other cofinality
(again, since $\kappa\ge \aleph_2$ there is such $\delta$ with
$\cf(\delta) \in \{\aleph_0,\aleph_1\} \backslash \{\emptyset\})$.
Now $C_{\mlt}$ is clearly a Dedekind cut of $K_\kappa$ hence Lemma
\ref{cutsymmetry} applies to it, too, but its first conclusion fails
for $\alpha = \delta$ hence its assumption fails.
So for some $\beta < \kappa$, the cut $(C_{\mlt}) \rest 
K_\beta$ of $K_\beta$ induces $C_{\mlt}$.  So for some 
increasing sequence $\langle \alpha_\varepsilon:\varepsilon < 
\kappa\rangle$ of ordinals $< \kappa$ and $c^+_\varepsilon \in K_\beta$ we have
 $b^+_{\alpha_{\varepsilon +1}}/b^-_{\alpha_{\varepsilon +1}} <
c^+_\varepsilon \le
b^+_{\alpha_\varepsilon}/b^-_{\alpha_\varepsilon}$ so $\langle
 c^+_\varepsilon:\varepsilon < \kappa\rangle$ is decreasing unbounded
 from below in $(C_{\mlt})^+ \cap K_\beta$, so unbounded in
 $(C_{\mlt})^+$.

This exemplifies $\kappa < d(K_\beta)$ but by \ref{3.4} we have
$d(K_\beta) \le h^+(K_\beta)$ but clearly $h(K_\beta) \le h(K)+ \beta$
hence $|h(K_\beta)| \le |h(K) + (\beta) < h^+(K) + \kappa = \kappa$, a
contradiction. 
\medskip

\noindent
\underline{Case IV: $C$ is a positive Dedekind cut, but not a Scott cut}

One argues as in the preceding case, considering $C_\add$
and using Lemma \ref{dedekind},
which leads to a symmetric additive cut and thus a contradiction
to the previous case.
In details choose $\langle b^-_\alpha:\alpha<\kappa\rangle,
\langle b^+_\alpha:\alpha<\kappa\rangle$ as in case III,
so $\langle b^+_\alpha-b^-_\alpha:\alpha<\kappa\rangle$ is a
\dec\ \seq\ in $C^+_{\rm add}$ unbounded from below in it.

If the cofinality of $C_{\rm add}$ from below is also $\kappa$,
recall that $C_{\rm add}$ is an additive cut by \ref{2.3.1}(2)
contradiction to case III.  If not, we repeat the argument in the end
of Case III.
\medskip

\noindent
\underline{Case V: $C$ is a (Dedekind not Scott) cut of $K_\kappa$}

Choose $a \in C^+$ hence $(a + C^-,a + C^+)$ is a positive cut of
$K_\kappa$ so we get a contradiction by Case IV.

As no case remains, Proposition \ref{termination}(1) is proved, and
thus the construction of a symmetrically complete extension
terminates.

As for clause (i) of \ref{termination}(2),
to estimate the cardinality of the resulting symmetrically
complete extension, recall that it has height at most $h(\kappa) +
\kappa +1$ hence $|h(K_{\kappa +1})| \le 
\kappa'=\max(h^+(K),\aleph_2)\le \max(|K|^+,\aleph_2)$  and hence
$K_{\kappa +1}$ has cardinality at most $2^{\kappa'}$. 
Moreover, similarly for any $\alpha<\kappa',|K_\alpha| \le 
2^{h^+(K)+\aleph_1}$ hence

\par \noindent
$|K_\kappa|=|\bigcup\limits_{\alpha<\kappa} K_\alpha|\leq
\sum\limits_{\alpha<\kappa} |K_\alpha| \leq
\sum\limits_{\alpha<\kappa} 2^{h^+(K)+\aleph_1}=\kappa'
+2^{h^+(K)+\aleph_1}=2^{h^+(K)+\aleph_1}$.  As $K_\kappa$ is dense in
$K_{\kappa +1}$ and $d(K_\kappa) \le h^+(K)$ also $|K_{\kappa +1}| \le
2^{h^+(K) + \aleph_1}$.

For clause (ii) of \ref{termination}, we define an
embedding $h_\alpha$ of $
K_\alpha$ into $K'$, increasing continuous
with $\alpha$ for $\alpha\leq \kappa$. For $\alpha=0,
h_0$ is the identity, for $\alpha$ limit take union and
for $\alpha=\beta+1$ use \ref{symmetry}.

Clauses (iii),(iv),(v) of \ref{termination} is easy too.
%hfill \qedref{termination}

\begin{discussion}
\label{d15}
How do we prove that the bound in \ref{termination} is right?  It goes
as in \cite[\S2]{Sh:405} but using a given decreasing sequence of
length $\theta$ for $\theta < \kappa$.  This is to be filled.
\end{discussion}
\newpage
 
\section{Concluding remarks}

\begin{discussion}
\label{e5}
It should be clear that there are considerably more general types of
closure that can be constructed in a similar manner. Let $\Theta$
be a class of possible cofinalities of cuts, that is pairs of regular
cardinals, and suppose that $\Theta$ is symmetric in the sense that
$(\theta_1,\theta_2)\in \Theta$ implies
$(\theta_2,\theta_1)\in \Theta$. Then we may consider
$\Theta$-constructions in which maximal independent sets of cuts,
all of whose cofinalities are restricted to lie in $\Theta$,
are taken. In order to get such a construction always to terminate, all that
is needed is the following:
\mn
\begin{enumerate}
\item[$(a)$] for all regular $\theta_1$, there is $\theta_2$
such that the pair $(\theta_1,\theta_2)$ is not in $\Theta$.
\end{enumerate}
\mn
From this it follows that:
\mn
\begin{enumerate}
\item[$(b)$]  for some regular $\kappa \geq h^+(K)+\aleph_2$, for every
$\theta_1$ regular $<\kappa$ there is a regular $\theta_2<\kappa$
\st\ $(\theta_1,\theta_2)\notin \Theta$.
\end{enumerate}
\mn
The proof is as above; in the symmetric case, $\Theta_{\rm sym}$
consists of all pairs $(\theta,\theta)$ of equal regular cardinals.
Clearly we may have to make the closure as large as we need $\kappa$ as in 
(b) above. Also in the proof of
\ref{termination}(1) in the multiplicative case we choose $\delta < \kappa$
\st\ $(\aleph_0,\cf(\delta))\notin \Theta$ and in the additive case,
we choose $\delta < \kappa$ such that $(\theta,\cf(\delta)) \notin \Theta$
(but of course change the cardinality bound).

Under the preceding mild conditions, such a $\Theta$-construction
provides an ``atomic'' extension of the desired type.
So we have $\Theta$-closure, and it is prime (as in
clause (ii) of \ref{termination}(2)).
We also can change the cofinality of $K$.  Below we elaborate.
\end{discussion}

\begin{definition}
\label{5.1}
Assume $\Theta$ is a set or class of pairs of regular
cardinals which is symmetric i.e.,
$(\kappa,\lambda)\in \Theta \Rightarrow (\lambda,\kappa) \in \Theta$.

\noindent
0) A $\Theta$-cut of $K$ is a cut $C$ with $\cf(C)\in \Theta$.

\noindent
1) We say that a real closed field is $\Theta$-complete \Iff \, 
no Dedekind cut of $K$ has cofinality from $\Theta$.

\noindent
2) We say that $L$ is a $\Theta$-hull of $K$ \when \, 
there is a maximal subset ${\cC}$ of $\cut_{\Theta}(K)
:= \{C:C$ a Dedekind cut of $K$
with cofinality $\in \Theta\}$ which is independent and
$a_C \in L$ for $C\in {\cC}$ realizing $C$ \st\ $L$
is the real closure of $K\cup \{a_C: C\in {\cC}\}$.

\noindent
2A) We say that $L$ is a weak $\Theta$-hull of $K$ \when \, there 
is a subset ${\cC}$ of ${\rm cut}_\Theta (K)$ which is 
independent and $a_C\in L$ for $C \in {\cC}$ realizing 
$C$ \st\ $L$ is the real closure of $K\cup \{a_C:C\in {\cC}\}$.

\noindent
3) We say that $\langle K_\alpha:\alpha \le \alpha^*\rangle$ is 
an associated $\Theta-\alpha^*$-chain over $K$ \when \,:
$K_\alpha$ is a real closed field, \incr\ \cont\ with $\alpha,
K_0=K$ and $K_{\alpha+1}$ is a $\Theta$-hull of $K_\alpha$
for $\alpha<\alpha^*$.

\noindent
3A) We say that $\langle K_\alpha:\alpha\leq \alpha^*\rangle$
is a weak associated $\Theta-\alpha^*$-chain over $K$ \when \,:
$K_\alpha$ is an increasing continuous sequence of real closed fields,
$K=K_0$ and $K_{\alpha+1}$ a weak $\Theta$-hull of $K_\alpha$
for $\alpha<\alpha^*$. We may omit $\Theta$ meaning
$\{(\kappa,\lambda):\kappa,\lambda$ are regular infinite
cardinals$\}$. We may omit ``over $K$".

\noindent
4)  We say that $\langle K_\alpha:\alpha\leq \alpha^*\rangle$ is an 
associated $\Theta$-chain over $K$ \when \,:
it is an associated $\Theta-\alpha^*$-chain over $K$,
$\alpha<\alpha^*\Rightarrow K_{\alpha+1} \ne K_\alpha$
and $K_{\alpha^*}$ is $\Theta$-complete.

\noindent
5)  Let $d'(K)$ be the minimal regular cardinal $\kappa$ (so infinite) \st:
for every non-Scott Dedekind cut $C$ of $K$ both cofinalities
of $C$ are $<\kappa$).
\end{definition}

\begin{theorem}
\label{5.2}
Let $K$ be a real closed field and $\Theta$ be as in
Definition \ref{5.1}.

\noindent
1) There is a $\Theta$-hull $L$ of $K$, see Definition \ref{5.1}(1).

\noindent
2) $L$ in (1) is unique up to \isom\ over $K$, and $K$ is cofinal in it.

\noindent
3) For every ordinal $\alpha^*$ there is an associated 
$\Theta-\alpha^*$-chain over $K$, see Definition \ref{5.1}(2)
and it is unique up to \isom\ over $K$.

\noindent
3A)  If $\langle K_\alpha:\alpha\leq\alpha^*\rangle$ is a weak 
associated $\Theta-\alpha^*$-chain over $K$ \then \, $K$ is cofinal in 
all its members in particular in $K_{\alpha^*}$.

\noindent
4)  If $\theta=\cf(\theta)<d(K)$ moreover $K$ has 
a non-Scott Dedkind cut of lower cofinality $\theta$ and
$(\forall \lambda) [(\theta,\lambda)\in \Theta]$
\then \, there is no associated $\Theta$-chain over $K$,
moreover no $\Theta$-complete extension of $K$.

\noindent
5) $\aleph_0 \le d' (K)\leq d (K)$.

\noindent
6)  If $d'(K)=\aleph_0$ \then\ $K$ is isomorphic to a sub-field of the
 field of reals.

\noindent
7)  If $\langle K_\alpha:\alpha\leq \alpha^*\rangle$ is a 
weak associated $\Theta-\alpha^*$-chain over $K$ \then \, $d'(K) \le 
\Max\{|\alpha^*|^+, d'(K)\}$.  

\noindent
7A)  In (7) also $|K_{\alpha^*}|\leq |K|+
|\alpha^*|+\sum \{|K|^\theta:(\theta,\theta)\in \theta$ and
$\theta < d'(K)\}$.

\noindent
8)  Assume $\kappa$ satisfies ``$\kappa \ne \cf(K),\kappa =
\cf(\kappa)>\aleph_0$ and for every regular $\theta<d'(K)+\kappa$ 
for some regular $\lambda<\kappa$ we have $(\theta,\lambda)\notin \Theta$,
\then\ there is an associated $\Theta$-chain over $K$ of
length $\leq \kappa$ (compare with part (4)).

\noindent
9)  If ($\forall$ regular $\theta$) ($\exists$ regular $\lambda$)
[$(\theta,\lambda)\notin \Theta$] \then\ for every
real closed field $K'$ there is an associated $\Theta$-chain over it.

\noindent
10)  If there is an associated $\Theta$-chain
$\langle K_\alpha:\alpha\leq \alpha^*\rangle$ over $K$ \then\:
\mn
\begin{enumerate}
\item[$(a)$]  $K_{\alpha^*}$ is a real closed field,
$\Theta$-complete extending $K$
\sn
\item[$(b)$] [universally]] if $K'$ is a $\Theta$-complete real closed field
extending $K$ \then\ $K_{\alpha^*}$ can be embedded into $K'$ over $K$
\sn
\item[$(c)$] [uniqueness]  if $L$ is a $\Theta$-complete real closed
field extending $K$ which can be embeded over $K$ into $L'$ for 
every $\Theta$-complete, $\Theta$-complete real closed field 
extending $K$ \then\ $L$ is \isom\ to $K_{\alpha^*}$ over $K$
\sn
\item[$(d)$]  in Clause (c), $|L| \le 2^{<\kappa}$ when 
$\kappa$ minimal as in (8).
\end{enumerate}
\end{theorem}

\begin{PROOF}{\ref{5.2}}
1)  As in \ref{b36}(1).

\noindent
2)  As in \ref{b36}(1).

\noindent
3)  Follows form (1)+(2).

\noindent
3A)  Easily by induction on $\alpha^*$.

\noindent
4) Let $C$ be a non-Scott cut of $K$ \st\ $\cf(C) = 
(\theta,\lambda)$, \wilog \, $C$ is positive and $\langle a_\alpha:
\alpha<\theta\rangle$ be an \incr\ \seq\ of members of $C \cap K_+$
cofinal in it.
As $C$ is non-Scott clearly for some $c\in K_+$
we have $a\in C^-\Rightarrow a+c\in C^-$ hence \wilog\
\mn
\begin{enumerate}
\item[$(*)$]   $c \in K_+$ and $a_\alpha+c<a_{\alpha+1}$ for 
$\alpha<\theta$.
\end{enumerate}
\mn
Now if $L$ is a $\Theta$-complete extension of $K$ then
there is a cut $C_1$ of $L$ \st\ $C^-_1= \{a\in L:a<a_\alpha$ for 
some $\alpha<\theta\}$.
Let $\cf(C_1)=(\lambda_1,\lambda_2)$, so necessarily
$\lambda_1=\theta$ and obviously $C^+_1$ cannot have
a first element $b$ as then $b-c\in C^+_1$ by (*),
so $\lambda_2\geq \aleph_0$. So by an assumption
$(\lambda_1,\lambda_2)\in \Theta, \text{ contradiction to } 
``L \text{ is } \Theta-\text{complete}"$.

\noindent
5)  Trivial, see Definition \ref{5.1}(5).

\noindent
6)  Easily $K$ is complete hence it is well known to be isomorphic to the 
field of reals.

\noindent
7)  So
\mn
\begin{enumerate}
\item[$(*)$]  $(a) \quad \langle K_\alpha:\alpha\leq \alpha^* \rangle$
is a weak $\Theta$-\seq
\sn
\item[${{}}$]  $(b) \quad K_0$ is unbounded in $K_{\alpha_*}$.
\end{enumerate}
\mn
Now we repeat the proof of \ref{3.4}.  The only place we have to say
more is why $\langle B^\varepsilon_i:i < \kappa,f(i) > \kappa\rangle$
has a constant end segment.  Otherwise, by the induction hypothesis,
its limit is a Scott cut of $L_\varepsilon$, but then $\langle a_i:i <
\kappa\rangle$ is cofinal if $B^-,B^-$ a Scott cut of $K_\kappa$.

\noindent
7A) Recall that $\cut_{\{(\theta,\lambda):\theta \ne 
\lambda {\rm regular}\}} (K)$ has cardinality $\leq |K|$.

\noindent
8) So assume $\langle K_\alpha:\alpha \le \kappa +1 \rangle$ is an
   associated $\Theta-(\kappa +1)$-chain over $K$.

First
\mn
\begin{enumerate}
\item[$(*)_1$]   like \ref{induced}, \ref{cutsymmetry} 
replacing ``symmetric'' cut by ``$\Theta$-cut''.
\end{enumerate}
\mn
Let
\mn
\begin{enumerate}
\item[$\boxplus$]   $C$ be a cut of $K_\kappa$ which is a 
$\Theta$-cut, we fix it for awhile.
\end{enumerate}
\mn
If for some $\alpha < \kappa$ the $\cat(C_\alpha) = C \rest K_\alpha$
induce $C$ on $K_\kappa$, then $C\rest K_\alpha$
is a cut of $K_\alpha$ of cofinality the same as $C$, hence
a $\Theta$-cut of $K_\alpha$ hence is realized in
$K_{\alpha+1}$ by the construction, say by $a$,
contradiction to ``$C\rest K_\alpha$ induce $C$ on $K_\kappa$".  

Hence
\mn
\begin{enumerate}
\item[$(*)_2$]  for no $\alpha<\kappa$ does $C\rest K_\alpha$ induce 
$C$ on $K_\kappa$ (so the assumption of \ref{cutsymmetry} holds). 
\end{enumerate}
\mn
This means that for every $\alpha<\kappa$ some $a_\alpha\in K
\setminus K_\alpha$ realizes $C$, so $a_\alpha \in C^-\vee a_\alpha
\in C^+$ so as we can replace $C$ by $|\{-b:b\in C^+\}, \{-b:b\in C^-\}|$
\wilog\ for arbitrarily large $\alpha<\kappa,a_\alpha \in C^-$, so
\mn
\begin{enumerate}
\item[$(*)_3$]  $a_\alpha \in C^-\subseteq K_\kappa$
realizes the cut $C\rest K_\alpha$ 
\sn
\item[$(*)_4$] $\cf(C^-)=\kappa$ so let $\cf(C)=(\kappa,\lambda)$.
\end{enumerate}
\medskip

\noindent
\underline{Case 1}:
$\lambda \ne \kappa$

Let $\sigma$ be regular $<\kappa$ \st\ $(\sigma,\lambda)\notin \Theta$.
Let $\{b_\beta:\beta<\lambda\rangle$ be a \dec\ \seq\
in $C^+$ unbounded from below in it so \wilog \,
\mn
\begin{enumerate}
\item[$(*)$]  $\{b_\beta:\beta<\lambda\} \subseteq K_{\alpha(*)}$.
\end{enumerate}
\mn
Now for some club $E$ of $\kappa$ we have
\mn
\begin{enumerate}
\item[$(*)$]  if $\delta\in E$ then $\delta>\alpha
(*)$ and $\{a_\alpha:\alpha<\delta\}$ is an unbounded subset of 
$C^-\cap K_\delta$ hence $\cf(C\rest K_\delta)=(\cf(\delta),\lambda)$.
\end{enumerate}
\mn
Choose $\delta\in E$, \st\ $\cf(\delta)=\sigma$.
But then $a_\delta$ realizes a $\{(\cf(\delta),\lambda)\}$-cut i.e.,
$\{(\sigma,\lambda)\}$-cut i.e., $\{ (\sigma,\lambda)\}$-cut which is a non
$\Theta$-cut by the choice of $\sigma$, \contr\ to $(*)_1$.
\medskip

\noindent
\underline{Case 2}:\ $\lambda=\kappa$ 
\newline
We repeat the proof of \ref{termination} after
\ref{cutsymmetry}, in (Case I) using $\kappa \neq \cf(K)$.

\noindent
9) For the given field $K'$ define $\kappa_n$ by induction on $n<\omega$

\[
\kappa_0 =|K'| \ (\hbox{or } d'(K))
\]

\mn
$\kappa_{n+1} = \Min\{\kappa:\kappa$ regular and if
$\theta<\kappa_n$ then for some $\lambda < \kappa$ we have 
$(\theta,\lambda) \notin \Theta\}$.

Now $K', (\sum \{\kappa_n:n<\omega\})^+$ satisfies the condition in (8).

\noindent
10)  Should be clear.
\end{PROOF}

\noindent
What about $\cf(K)$, we have not changed it in all our completion. 
It doesn't make much difference because
\begin{claim}
\label{10.4}
1)  For $K$ and regular $\kappa (\ge \aleph_0)$ there is $L$ such that:
\mn
\begin{enumerate}
\item[$(a)$]  $K\subseteq L$
\sn
\item[$(b)$]  $\cf (L)=\kappa$
\sn
\item[$(c)$]  if $K\subseteq L'$ and $\cf(L')=\kappa$ \then \, we 
can embed $L$ into $L'$ over $K$
\sn
\item[$(d)$]  $(\alpha) \quad$ if $\cf(K) = \kappa$ then $K=L$
\sn
\item[${{}}$]  $(\beta) \quad$ if $\cf(K) \ne \kappa$ then in clause (c)
we can add: there is an embedding $f$

\hskip25pt  of $L$ into $L'$
over $K$ \st\ ${\rm Rang} (f)$ is unbounded in $L'$.
\end{enumerate}
\mn
2) We can combine this with \ref{5.2}.
\end{claim}

\begin{PROOF}{\ref{10.4}}
Should be clear.
\end{PROOF}
\bigskip

\centerline {$* \qquad * \qquad *$}
\bigskip

\noindent
We have concentrated on real closed fields.  This is justified by
\begin{claim}
\label{5.6}
1)  Assume $K$ is an ordered field, $\theta=\cf(K)$ and $K$ has no
$\{(\theta,\theta)\}$-cut. \Then\ $K$ is real closed.

\noindent
2) In Theorem \ref{5.2} if we add $(\cf(K),\cf(K))\notin \Theta$ and
deal with ordered fields, it still holds.
\end{claim}

\begin{claim}
\label{5.7}
1)  If $F$ is an ordered field of cardinality $\mu > \aleph_0$ 
\then \, there is $F'$ \st\ :
\mn
\begin{enumerate}
\item[$(a)$]  $F'$ is a real closed field of cardinality $\lambda$
\sn
\item[$(b)$]  $F'$ extends $F$
\sn
\item[$(c)$]  if $C$ is a Dedekind cut of $F'$ of
cofinality $(\theta_1,\theta_2)$ then $\theta_1=\theta_2$
\sn
\item[$(d)$]  if $F''$ is another real satisfying (a),(b),(c) then
$F'$ can be embedded into $F'$ over $F$.
\end{enumerate}
\end{claim}

\begin{PROOF}{\ref{5.7}}
By \ref{5.2} applied to $\Theta = \{(\theta_1,\theta_1):\theta_1 \ne
\theta_2$ are regular (so infinite)$\}$.
\end{PROOF}

\noindent
Observe
\begin{claim}
\label{b13}
Assume $\langle K_\alpha:\alpha \le \gamma\rangle$ is increasing and
$C$ a cut of $K_\gamma$ and let $C_\alpha = C \rest K_\alpha$ for
$\alpha \le \gamma$.

\noindent
1) If $1 \in C^-$ (or just $C^- \cap (K_\alpha)_+ \ne \emptyset$ \then \,
   $C_{\add} \rest K_\alpha = (C_\alpha)_{\add}$.

\noindent
2) If $C$ is an additive cut \then \, each $C_\alpha$ is an additive
   cut and $C_{\add} \rest K_\alpha = (C_\alpha)_{\add}$.
\end{claim}

\begin{PROOF}{\ref{b13}}
Should be clear.
\end{PROOF}

\begin{claim}
\label{b15}
If $\bar K = \langle K_\alpha:\alpha \le \alpha_* \rangle$ is chains
over $K_0$ \then \, $h(K_{\alpha_*}) \le h(K_0) + \alpha_*$.
\end{claim}

\begin{PROOF}{\ref{b15}}
Implicite in \S3.
\end{PROOF}

\begin{remark}
Used in the end of \S4.
\end{remark}

\noindent
So (see \S0) we wonder
\begin{question}
\label{b22}
For $T$ dependent model which are $|T|^+$-saturated $\kappa$-saturated
for types which does not split over sets of cardinality $\le |T|$:
\mn
\begin{enumerate}
\item[$(a)$]  what \cite[\S5]{Sh:715} gives
\sn
\item[$(b)$]  when do we have symmetric cuts? (see
\cite[\S2]{Sh:405}).
\end{enumerate}
\end{question}

\begin{question}
\label{b24}
Similarly for o-minimal theory $T$.
\end{question}
\newpage

\section {Symmetric cuts exist for BPA}

We like in a model of, e.g. BPA to immitate \cite[\S2]{Sh:405}, see
\ref{y6} in \cite{Sh:757}.  The answer is essentially that it suffices
to have a function $f$ where $f(x)$ behaves like $x^x$.

We still have to sort out the beginning of the induction.  If we ask
only to have some countable model $M_*$ of $T$ such that ``if $M_* \prec
M$ \then \, $M$ has a symmetric cut", then this is O.K.

But it seems too much to ask for the existence of
such function.  So we may try to rework the proof using an almost
function: there is a formula giving a convex set of possible values,
and its non-existent infinum (i.e. exists in the completion) is as
required.  

\begin{theorem}
\label{f2}
The model $N$ has a symmetric cut when:
\mn
\begin{enumerate}
\item[$(a)$]  $(\alpha) \quad N$ is a model of $T$
\sn
\item[${{}}$]  $(\beta) \quad T$ is a first order theory 
extending the theory of ordered semi-rings 

\hskip25pt (so have 0,1 order, addition and multipliation
with the usual rules 

\hskip25pt but $x-y$ does not necessarily exist), which may
have additional 

\hskip25pt symbols (e.g. $x^x$)
\sn
\item[$(b)$]  for some formula $\varphi(x,y) \in \bbL(\tau_T)$ the
theory $T$ implies:
\sn
\item[${{}}$]  $(\alpha) \quad 0 < x \rightarrow (\exists y)\varphi(x,y)$
\sn
\item[${{}}$]  $(\beta) \quad \varphi(x,y) \rightarrow x < y$
\sn
\item[${{}}$]  $(\gamma) \quad 0 < x \rightarrow (\exists y)(\forall
  z)(\varphi(x,z) \leftrightarrow z \le y)$
\sn
\item[${{}}$]  $(\delta) \quad 0 < x_1 < x_2 \wedge \varphi(x_1,y_1) \wedge
\varphi(x_2,y_2) \rightarrow x_1 y_1 < y_2$
\sn
\item[$(c)$]  there are $e_n \in N$ for $n \in \bbN$ such that for
every $n,m \in \bbN$ we have: $\varphi(e_{n+1},x)
  \rightarrow x \le e_n$ and $m e_{n+1} < e_n$.
\end{enumerate}
\end{theorem}

\begin{remark}
1) If $T = \BPA$ then:
\mn
\begin{enumerate}
\item[$(a)$]  $\varphi(x,y) = (y = x^x)$ satisfies clause (b) of
  \ref{f2}
\sn
\item[$(b)$]  if $N \models T$ is non-standard, $T$ is as above, \then
  \, there are non-standard $e_n$ as in clause (c): choose $e_0$
  arbitrarily, then $e_{n+1}$ is the minimal $x$ such that $x^{(x^x)}
  \ge e_n$.
\end{enumerate}
\end{remark}

\begin{PROOF}{\ref{f2}}
So let $T,\varphi$ be as in the assumptions of \ref{f2}(2) and let $N$
be a model of $T$.  As satisfaction is all this proof is in $N$ we
omit $N \models$; and $a,b,c,d$ will denote members of $N$; we adopt the
notation:
\mn
\begin{enumerate}
\item[$\boxplus_1$]  $(a) \quad a <_\varphi b$ means $(\forall
  x)(\varphi(a,x) \rightarrow x \le b)$
\sn
\item[${{}}$]  $(b) \quad a \ll b$ means $n \cdot a < b$ for every
  $n$.
\end{enumerate}
\mn
By induction on the limit ordinal $\beta$ we will choose elements
$a_{\alpha,n},b_{\alpha,n}$ in $N$ for $n \in \bbN$ and $\alpha <
\beta$ such that:
\mn
\begin{enumerate} 
\item[$\boxplus^2_\beta$]  for all $n$ and for all $\alpha_1 < \alpha_2 <
\beta$ we have
\sn
\item[${{}}$]  $(a) \quad 1 < a_{\alpha_1,n} < a_{\alpha_2,n} 
< b_{\alpha_2,n} < b_{\alpha_1,n}$
\sn
\item[${{}}$]  $(b) \quad b_{\alpha_1,n+1} <_\varphi a_{\alpha_1 +1,n}
  -  a_{\alpha_1,n}$.
\end{enumerate}
\mn
\underline{Case 1}:  $\beta = \omega$

By clause (c) of the assumptions of Theorem \ref{f2} there are $e_n$
for $n < \omega$ as there, let $d_n = e_{2n+1},d^*_n = e_{2n}$, so
$d_n,d^*_n$ for $n < \omega$ are such that $d_{n+1} <_\varphi d^*_{n+1} \ll
d_n$ for $n$.

We let $a_{i,n} = d_{n+1} + i \cdot d^*_{n+1}$ and $b_{i,n} = 
d_n-i-1$.  We should check that $\boxplus^2_\beta$ holds.  Now clause
(a) holds because: $a_{i,n} < a_{i+1,n}$ as $d^*_{n+1} > 0$ and
$a_{i+1,n} - a_{i,n} = d^*_{n+1}$; and $a_{i,m} < b_{i,m}$ as $a_{i,m}
= d_{m+1} + id^*_{m+1} \le (i+1) d^*_{m+1} < (i+1+i) d^*_{m+1}-i-1 < 
d_m - i-1 = b_{i,m}$; and lastly $b_{i+1,n} < b_{i,n}$ as 
$d_{n+i+1} < \ldots < d_n$ hence $i+1 < d_n$.

Also clause (b) of $\boxplus^2_\beta$ holds $b_{i,n+1} = d_{n+1} -i-1
\le d_{n+1} <_\varphi d^*_{n+1} = a_{i+1,n} - a_{i,n}$.
\medskip

\noindent
\underline{Case 2}:  $\beta$ is a limit ordinal such that $\gamma <
\beta \Rightarrow \gamma + \omega < \beta$.

There is nothing to do.
\medskip

\noindent
\underline{Case 3}:  $\beta = \gamma + \omega,\beta$ a limit ordinal

So the elements
$a_{\alpha,n}$ and $b_{\alpha,n}$ have been chosen for $\alpha <
\gamma$ with $\gamma$ a limit ordinal.  Let $A_n,B_n$ be the ranges of
the sequences $a_{\alpha,n},b_{\alpha,n}$ (for $\alpha < \gamma$)
respectively so $A_n < B_n$, i.e. $(\forall a \in A_n)(\forall b \in
B_n)(a <_N b)$.  If one of the pairs $(A_n,B_n)$ determines a cut in
$N$, then it defines a cut as a desired (i.e. the cut
$(\{x:(\exists y \in A_n)(x<y)\},\{x:(\exists y \in B_n)(y < x\}$).

Assume therefore:
\mn
\begin{enumerate}
\item[$\boxplus_3$]  for all $n$ there is an element $c_n$ with $A_n <
c_n < B_n$, i.e. $(\forall x \in A_n)(\forall y \in B_n)[x < c_n < y]$.
\end{enumerate}
\mn
Under this assumption we will continue the construction by defining
$a_{\gamma+i,n}$ and $b_{\gamma +i,n}$ for all finite $i,n$ thus
finishing this case too (hence the proof).

We set
\mn
\begin{enumerate}
\item[$(*)_1$]  $(a) \quad$ let $b^*_{\alpha,n}$ be such that
  $\varphi(b_{\alpha,n},b^*_{\alpha,n})$ for $\alpha < \gamma,n <
  \omega$
\sn
\item[${{}}$]  $(b) \quad$ let $c^*_n$ be such that
  $\varphi(c_n-1,c^*_n)$
\sn
\item[${{}}$]  $(c) \quad$ let $c^{**}_n$ be such that 
$\varphi(c_n,c^{**}_n)$
\sn
\item[${{}}$]  $(d) \quad c'_n := c_n - c^{**}_{n+1}$.
\end{enumerate}
\mn
Now we observe that
\mn
\begin{enumerate}
\item[$(*)_2$]  $A_n < c'_n$ (i.e. $(\forall x \in A_n)(x < c'_n))$.
\end{enumerate}
\mn
[Why?  Since for every $n \in \bbN$ and $\alpha < \gamma$ we have:
$a_{\alpha,n} < a_{\alpha +1,n} - b^*_{\alpha,n+1} < a_{\alpha+1,n} - 
c^{**}_{n+1} < c_n - c^{**}_{n+1} = c'_n$.]

Why?  First inequality by $\boxplus^2_\gamma(b)$ recalling
  $\boxplus_1(a)$.  Second inequality
as $b_{\alpha,n+1} > c_{n+1}$ by $\boxplus_3$ and clause
$(b)(\delta)$ of the theorem (check).  Third inequality holds as
$a_{\alpha +1,n} < c_n$ by $\boxplus_3$.  Fourth equality by the
choice of $c'_n$.]

We set (for $i,n < \omega$):
\mn
\begin{enumerate}
\item[$(*)_3$]  $(a) \quad a_{\gamma +i,n} := c'_n + i \cdot c^*_{n+1}$
\sn
\item[${{}}$]  $(b) \quad b_{\gamma +i,n} := c'_n + c_{n+2} \cdot
  c^*_{n+1} - i$.
\end{enumerate}
\mn
So we have to check the inductive demands, this means
\medskip

\noindent
\underline{Clause (a) of $\boxplus^2_{\gamma+1}$}:  This means that
the following inequalities hold below $i,n \in \bbN$ and $\alpha < \gamma$:
\medskip

\noindent
\underline{$a_{\alpha,n} < a_{\gamma+i,n}$}:
\medskip

\noindent
Why?  As $a_{\alpha,n} < c'_n$ by $(*)_2$ and $c'_n < c'_n +
i \cdot (c^*_{n+1})$ because $c^*_{n+1} > 1$ and $c'_n + i \cdot c^*_{n+1} 
= a_{\gamma+i,n}$ by the choice of the latter, see $(*)_3(a)$.
\medskip

\noindent
\underline{$a_{\gamma + i,n} < a_{\gamma +i+1,n}$}:

Why?  Trivial by $(*)_3(a)$ as $c^*_{n+1} > 0$.
\medskip

\noindent
\underline{$a_{\gamma+i,n} < b_{\gamma+i,n}$}:

Why?  By $(*)_3$ this means $i \cdot c^*_{n+1} < c_{n+2} \cdot
c^*_{n+1} -i$ but $i <  c_{n+2}$ hence $i \cdot c^*_{n+1} <
(c_{n+2} -i) \cdot c^*_{n+1} = c_{n+2} \cdot c^*_{n+1} - i \cdot
c^*_{n+1} < c_{n+2} \cdot c^*_{n+1} -i$ is as required.
\medskip

\noindent
\underline{$b_{\gamma+i+1,n} < b_{\gamma +i,n}$}:

Why?  Trivial by $(*)_3(b)$ because $i+1 < c'_n$.
\medskip

\noindent
\underline{$b_{\gamma+i,n} < b_{\alpha,n}$}:

Why?  Note that $c_n < b_{\alpha,n}$ by $\boxplus_3$ as $b_{\alpha,n}
\in B_n$ hence it suffices to prove $b_{\gamma +i,n} < c_n$, which by
$(*)_3(b)$ means $c'_n + c_{n+2} \cdot c^*_{n+1}-i < c_n$.  By the choice of
$c'_n$ in $(*)_1(d)$ this means $c_{n+2} \cdot c^*_{n+1}-i 
< c^{**}_{n+1}$ and as $c_{n+2} < (c_{n+1}-1)$ it suffices to prove
$(c_{n+1}-1) \cdot c^*_{n+1} < c^{**}_{n+1}$ 
which holds by the choice of $c^*_{n+1},c^{**}_{n+1}$ in $(*)_2(b),(c)$ see 
Clause $(b)(\delta)$ of the claim.
\medskip

\noindent
\underline{Clause (b) of $\boxplus^2_{\gamma +n}$}:

We have to show $b_{\gamma +i,n+1} <_\varphi a_{\gamma +i+1,n} -
a_{\gamma +i,n}$.

So assume $\varphi(b_{\gamma +i,n+1},b'_{\gamma +i,n+1})$ holds and we
have to show $b'_{\gamma +i,n+1} \le a_{\gamma +i+1,n} - a_{\gamma
  +i,n}$.  By clause $(*)_3(a)$ this means $b'_{\gamma +i,n+1} <
c^*_{n+1}$.  By the choice of $b'_{\gamma +i,n+1}$
above and of $c^*_{n+1}$ in $(*)_1(b)$ it suffices to prove (that the
model $N$ satisfies) $\varphi(b_{\gamma +i,n+1},y_1) \wedge
\varphi(c_{n+1} -1,y_2) \rightarrow y_1 < y_2$.  

By clause $(b)(\delta)$ of the claim, it suffices to 
prove $b_{\gamma +i,n+1} < c_{n+1} - 1$ which by $(*)_3(b)$ means
$c'_n + c_{n+2} \cdot c^*_{n+2} - i < c_{n+1}$ which by $(*)_1(d)$
means $c_{n+2} \cdot c^*_{n+1} - i < c^{**}_{n+1}$ which is proved above by
$(b)(\delta)$ of the claim.
\end{PROOF}
\newpage

 % PRIVATE PART 1 

\bibliographystyle{amsalpha}
\bibliography{shlhetal}

\end{document}